\documentclass{amsart}

\newtheorem{theorem}{Theorem}[section]
\newtheorem{proposition}[theorem]{Proposition}

\newtheorem{corollary}[theorem]{Corollary}

\newtheorem{definition}[theorem]{Definition}
\newtheorem{example}[theorem]{Example}

\newtheorem{remark}[theorem]{Remark}

\numberwithin{equation}{section}

\begin{document}

\title[Regularity bounds for curves]
{Regularity bounds for curves by minimal generators and Hilbert function}

\author{Francesca Cioffi}
\address{Dipartimento di Matematica e Applicazioni, 
Universit\`a di Napoli ``Federico II", 80126 Napoli, Italy}
\email{francesca.cioffi@unina.it}
\thanks{The authors were supported in part by MURST and GNSAGA}

\author{Maria Grazia Marinari}
\address{Dipartimento di Matematica, Universit\`a di Genova,
16146 Genova, Italy}
\email{marinari@dima.unige.it}

\author{Luciana Ramella}
\address{Dipartimento di Matematica, Universit\`a di Genova,
16146 Genova, Italy}
\email{ramella@dima.unige.it}


\date{}


\begin{abstract}
Let $\rho_C$ be the regularity of the Hilbert function of a projective 
curve $C$ in $\mathbb P^n_K$ over an algebraically closed field $K$ 
and $\alpha_1,\ldots,\alpha_{n-1}$ be minimal degrees for which there 
exists a complete intersection of type $(\alpha_1,\ldots,\alpha_{n-1})$ 
containing the curve $C$.
Then the Castelnuovo-Mumford regularity of $C$ is
upper bounded by $\max\{\rho_C+1,\alpha_1+\ldots+\alpha_{n-1}-(n-2)\}$. 
We study and, for space curves, refine the above bound providing several 
examples.
\end{abstract}

\maketitle

\section{Introduction}

Estimates for the Castelnuovo-Mumford regularity of a standard graded 
algebra $A$ ``give a measure of size of $A$" and it is interesting to 
predict such estimates from the Hilbert function of $A$ (\cite{V}, 
Remark 2.6). In this context, it is well known that the Hilbert polynomial 
of a graded polynomial ideal gives a bound for the Castelnuovo-Mumford 
regularity (Gotzmann's Regularity Theorem), but when the ideal is saturated 
this bound can be very far from being sharp. 

Here we tackle the question of bounding the Castelnuovo-Mumford regularity
of the algebra $A=S/I$, where $S=K[x_0,\ldots,x_n]$ and $I$ is the defining 
ideal of a curve $C$ over an algebraically closed field $K$ with general 
hyperplane section $Z$. 
We note that, by well known results about zero-dimensional schemes, an upper 
bound for the Castelnuovo-Mumford regularity $reg(C)$ of $C$ arises  
in terms of the regularity of the Hilbert function of $S/I$ and of the 
regularity of the Hilbert function of $Z$ (see Theorem 3.4). 
Then, involving also degrees of minimal generators of $I$, we study upper 
bounds for $reg(C)$, by means of upper bounds for $reg(Z)$ when $C$ is 
known and $Z$ is not (Proposition 4.1). 

Bounds in terms of the degrees of defining equations of
projective schemes are given for smooth \cite{BEL} and locally complete
intersection \cite{CU} schemes.

To address the above question, let $\mathbb P^n_K$ be a projective space of 
dimension $n$ over $K$ and $C$ be a curve 
in $\mathbb P^n_K$, i.e. a non-degenerate projective scheme of dimension $1$. 
It is well known that, denoting $\rho_C$ the regularity of the Hilbert 
function of $C$, $\omega_I$ the maximal degree of minimal generators of 
the defining ideal $I$ of $C$ and $reg(C)$ the Castelnuovo-Mumford regularity 
of $C$, then 
\begin{equation}\rho_C+1 \leq reg(C), \quad \omega_I \leq reg(C).
\end{equation} 

In this paper we note how it is possible to ``reverse" the inequalities (1.1)
for curves basing on 
Theorem 3.4 and results of \cite{DGO}, thus yielding that, if 
$\alpha_1,\ldots,\alpha_{n-1}$ are minimal degrees for which there 
exists a complete intersection (c.i. for short) of type $(\alpha_1,\ldots,\alpha_{n-1})$ 
containing the curve $C$, then (Proposition 4.1)
\begin{equation}
reg(C)\leq \max\{\rho_C+1,\alpha_1+\ldots +\alpha_{n-1}-(n-1)+1\}.
\end{equation}

Degrees of generators of the defining ideal, together with other
invariants, are used for bounding regularity of space curves 
in \cite{MMN,MM-R} and, in this context, it is 
worthwhile to see also \cite{N2} (see the references therein too). 

The paper is organized as it follows. 
In section 2 we set notation and recall some definitions and known
results on the Castelnuovo-Mumford regularity of a homogeneous ideal.
In section 3 we recall and state (Theorem 3.4) some information about 
the Castelnuovo-Mumford regularity of projective schemes that follows 
by the study of general hyperplane sections. 
Then, in section 4, by Proposition 4.1 
we note that formula (1.2) raises from well known results 
about zero-dimensional schemes of \cite{DGO}. Hence, we focus our attention on 
space curves, refining the above bound 
(Corollary 4.3) and comparing it with other bounds of the same type, 
as a bound for integral curves which is described by means of the degree 
of the curves and based on the shape of Borel ideals (Proposition 4.5).

Many of the computations have been performed by {\tt Points} \cite{Points} 
and \emph{{{\hbox{\rm C\kern-.13em o\kern-.07em C\kern-.13em o\kern-.15em A}}}}
\cite{CoCoASystem}.

\section{General setting}

Let $K$ be an algebrically closed field, 
$S=K[x_0,\ldots,x_n]$ the polynomial ring over $K$ in $n+1$ variables 
endowed with the deg-rev-lex term-orderig such that $x_0>x_1>\ldots>x_n$, 
$\mathbb P^n_K = Proj\  S$ the projective space of dimension $n$ over $K$.

If $I\subset S$ is a homogeneous ideal, $I_t$ denotes the subset of $I$ 
consisting of the homogeneous polynomials of $I$ of degree $t$ and 
$I_{\leq t}$ the subset of $I$ consisting of the homogeneous polynomials 
of $I$ of degree $\leq t$. Moreover, $\alpha_I$ is the {\it initial degree} 
of $I$ and $\omega_I$ is the {\it maximal degree} of minimal generators of $I$.

We will use freely the common notation of sheaf cohomology 
in ${\mathbb P}^n_K$, referring to \cite{Gr,Mi} for notation and 
basic results. 
If $\mathcal F$ is a coherent sheaf on $\mathbb P^n_K$ we will write 
$H^i(\mathcal F)$ instead of $H^i(\mathbb P^n_K,\mathcal F)$.

\begin{definition} \rm
A finitely generated graded $S$-module $M$
is $m$-{\it regular} if the $i$-th syzygy module of $M$ is generated in
degree $\leq m+i$, for all $i\geq 0$. The {\it regularity} $reg(M)$ of $M$ 
is the smallest integer $m$ for which $M$ is $m$-regular.
\end{definition}

\begin{remark}\rm
The above definition can be applied to a homogeneous ideal $I$ of $S$.
So, if $\omega_I$ is the maximal degree of minimal generators of $I$, 
then $\omega_I \leq reg(I)$.
\end{remark}

\begin{definition}\rm
The {\it saturation} of a homogeneous ideal $I$ is $I^{sat}=\{f\in S\ / \ \ 
\forall \ j=0,\ldots,n,\exists \ r \in \mathbb N : x_j^r f \in I\}$.
The ideal $I$ is {\it saturated} if $I^{sat}=I$. The ideal $I$ is 
{\it $m$-saturated} if $(I^{sat})_t = I_t$ 
for each $t\geq m$. The {\em satiety} 
of $I$ is $sat(I)=\min\{t \ / \ \text{$I$ is $t$-saturated}\}$.
\end{definition}

\begin{definition} \rm A coherent sheaf $\mathcal F$ on $\mathbb P^n$ is
$m$-{\it regular} if $H^i(\mathcal F(m-i)) = 0$ for all
$i>0$. The {\it Castelnuovo-Mumford regularity } (or {\it regularity})
of $\mathcal F$ is $reg(\mathcal F) = \min\{ m \ / \  \mathcal F$ is
$m$-regular $\}$.
\end{definition}

\begin{proposition} 
{\rm (\cite{Gr} Proposition 2.6)}
A homogeneous ideal $I$ is $m$-regular if and only if $I$ is $m$-saturated
and its sheafification $\mathcal I$ is $m$-regular. Hence, for a saturated
homogeneous ideal $I$, the regularity of $I$ equals the regularity of 
its sheafification.
\end{proposition}

\begin{definition}\rm
If $X \subset \mathbb P^n_K$ is a closed subscheme, its {\em regularity}
$reg(X)$ is the regularity $reg(I)$ of its defining ideal $I$.
\end{definition}

For a finitely generated graded $S$-module $M$ we let $H_M(t):=dim_K M_t$ be
the {\it Hilbert function} of $M$, $\Delta H_M(t) := H_M(t) - H_M(t-1)$, 
for $t\geq 1$, and $\Delta H_M(0):=1$. 

We recall that, for $t>>0$, $H_M(t)=P_M(t)$ 
where $P_M(z) \in K[z]$ is the {\it Hilbert polynomial} of $M$.
The {\it regularity} of $H_M$ is $\rho_M := \min\{\bar t \ / \ H_M(t)=P_M(t), 
\forall t\geq \bar t\}$.

If $I$ is the defining ideal of a closed subscheme $X\subset \mathbb P^n_K$ of
dimension $k$, instead of $H_{S/I}$, $P_{S/I}$, $\rho_{S/I}$ 
we can write $H_{X}$, $P_{X}$, $\rho_X$. Recall that the 
Hilbert polynomial of $X$ is also 
$P_X(t) = \sum _{i=0}^{k}(-1)^i h^i({\mathcal O}_X(t))$ 
(see \cite{H} Exercise III.5.2).

\begin{remark}\rm 
Let $X\subset \mathbb P^n_K$ be a non-degenerate scheme of dimension $k$. 
From the short exact sequence $0\to {\mathcal I}_X(t) 
\to {\mathcal O}_{{\mathbb P}^n}(t)\to {\mathcal O}_X(t)\to 0$ we get 
the long cohomology exact sequence
$$0 \rightarrow H^0({\mathcal I}_X(t)) \rightarrow H^0({\mathcal
O}_{{\mathbb P}^n}(t)) \rightarrow H^0({\mathcal O}_X(t))
\rightarrow H^1({\mathcal I}_X(t))\rightarrow 
H^1({\mathcal O}_{{\mathbb P}^n}(t)) \to $$
$$\rightarrow H^1({\mathcal O}_X(t)) \rightarrow H^2({\mathcal I}_X(t))
\rightarrow \cdots \rightarrow H^n({\mathcal O}_{{\mathbb P}^n}(t)) 
\rightarrow 0 $$
from which it follows that:
\begin{itemize}
\item[(1)] $H_X(t) = h^0({\mathcal O}_{\mathbb P^n}(t))-h^0({\mathcal I}_X(t))=
h^0({\mathcal O}_C(t))-h^1({\mathcal I}_X(t))$;
\item[(2)] $h^{i+1}(\mathcal I_X(t)) = h^i(\mathcal O_X(t))$ for 
$1\leq i\leq k$ and for $t \geq -n$,
$h^i(\mathcal I_X(t))=0$ for every $i\geq k+2$;
\item[(3)] $reg(X) \geq \rho_X+1$;
\item[(4)] when $k=1$ we have $h^2(\mathcal I_X(t)) = h^1(\mathcal I_X(t))$ 
for every $t \geq\rho_X$.
\end{itemize}
\end{remark}

\section{Castelnuovo-Mumford regularity and hyperplane sections}

Let $h$ be a {\em general} linear form for a homogeneous ideal $I\subset S$, 
i.e. $h\in S_1$ is a homogeneous polynomial of degree $1$ which is not 
a zero-divisor on $S/I^{sat}$.

\begin{proposition} {\rm (\cite{BS}, Lemma 1.10)} A homogeneous ideal 
$I$ is $m$-regular if and only if $I$ is $m$-saturated and $(I,h)$ is 
$m$-regular.
\end{proposition}

\begin{remark}\rm \label{zero} 
\rm From the above proposition it follows immediately that
$reg(I) = \max\{sat(I),reg((I,h))\}$. Thus, if $X$ is an arithmetically 
Cohen-Macaulay (for short aCM) scheme of dimension $k$, we have 
$reg(X)=\rho_X+k+1$.
\end{remark}

\begin{remark} \label{uno} \rm
By the short exact sequence
$$0 \longrightarrow (S/I)_{t-1} \buildrel{\cdot h} \over \longrightarrow
(S/I)_t \longrightarrow (S/J)_t\longrightarrow 0$$
it holds that $H_{S/J}(t) = \Delta H_{S/I}(t)$ and then
$\rho_{S/J}=\rho_X+1$.
\end{remark}

From now on, let $X\subset \mathbb P^n_K$ be a scheme of dimension $k>0$, $I$ 
its defining ideal and $J=(I,h)$. Moreover, let $Z\subset\mathbb P^{n-1}_K$
be the scheme of dimension $k-1$ defined by the saturated ideal 
$J^{sat} = (I,h)^{sat}$, i.e. the general hyperplane section of $X$.

Note that, for every $t\geq \max\{\rho_Z,\rho_{S/J}\}$, it is 
$H_{S/J}(t)=H_Z(t)$, hence
$$sat(J)\leq \max\{\rho_Z,\rho_{S/J}\}.$$

The following result is obtained by an easy computation involving Hilbert 
functions and the definition of saturated ideal.

\begin{theorem} With the above notation 
$$reg(X) = \max\{\rho_X+1,reg(Z)\}.$$
Namely,
\begin{itemize}
\item [(1)] $\rho_X+1\geq reg(Z)$ $\Rightarrow$ $sat(J)=reg(X)=\rho_X+1$,
\item [(2)] $\rho_X+1=reg(Z)-1$ $\Rightarrow$ $sat(J)\leq reg(Z)-1$ and 
$reg(X) = reg(Z)$,
\item [(3)] $\rho_X+1<reg(Z)-1$ $\Rightarrow$ $sat(J)=reg(Z)-1$ and 
$reg(X) = reg(Z)$.
\end{itemize}
\end{theorem}

\begin{proof}
By Propositions 2.5 and 3.1 it follows that $reg(X)=
reg(I)=reg(J)=\max\{sat(J),reg(J^{sat})\}$. Hence,
$reg(X)=\max\{sat(J),reg(Z)\}$.

If $\rho_X+1\geq reg(Z)$ then $\rho_{S/J}=\rho_X+1\geq \rho_Z+1$, 
by Remarks 2.7(3) and 3.3. Hence, for every $t\geq \rho_{S/J}\geq\rho_Z+1$ 
we have $H_{S/J}(t)=P_Z(t)=H_Z(t)$ but $H_{S/J}(\rho_X)\not=P_Z(\rho_X)
=H_Z(\rho_X)$. Since $J\subseteq J^{sat}$, it follows that, for every 
$t\geq \rho_X+1$, $J_t=J^{sat}_t$, meanwhile $J_{\rho_X}\not=J^{sat}_{\rho_X}$. 
Thus $sat(J)=\rho_X+1$. 

If $\rho_X+1=reg(Z)-1$, with analogous arguments as above it holds $sat(J)\leq
\rho_X+1=\rho_Z$ and so $reg(X)=reg(Z)$.

If $\rho_X+1<reg(Z)-1$, similarly we have $H_{S/J}(\rho_Z-1)=P_Z(\rho_Z-1)$ 
while $H_Z(\rho_Z-1)\not=P_Z(\rho_Z-1)$ and $H_{S/J}(\rho_Z)=P_Z(\rho_Z)=
H_Z(\rho_Z)$, so $sat(J)=\rho_Z$. 
\end{proof}

Note that, by applying Theorem 3.4 to a curve $C\subset\mathbb P^n_K$ 
with general hyperplane section $Z$, we get 
$reg(C)=\max\{\rho_C+1,\rho_Z+1\}=\min\{t\geq \rho_Z+1 \ \vert \ 
\Delta H_C(t)=0\}$, being $reg(Z)=\rho_Z+1$ by Remark 3.2. So, bounds for 
zero-dimensional schemes' regularity enter into our aims in section 4.

\begin{remark} \label{due} \rm In \cite{Ba} it is proved that 
if $C$ is an integral curve of degree $deg(C)$ then
$\rho_Z \leq \Bigl\lceil \frac{deg(C)-1}{n-1} \Bigr\rceil$.
\end{remark}

\begin{example} \rm We give examples of curves for every case of Theorem 3.4.
In (2) and (3) we apply the technique used in \cite{MMN} to construct 
curves with ``high degree generators", involving basic double linkages. 
We recall the procedure to obtain a basic double link of type 
$(a_1,\ldots,a_{n-1})$ from a curve $C$ of $\mathbb P_K^n$. \hskip 1mm 
Consider homogeneous polynomials $F_1,\ldots,F_{n-1}$ of $K[x_0,\ldots,x_n]$ 
of degrees $a_1,\ldots,a_{n-1}$ respectively, with $F_1$ general, 
$F_2,\ldots,F_{n-1} \in I$ and $(F_1,\ldots,F_{n-1})$ a regular sequence,
then the curve $C'\subset\mathbb P_K^n$ defined by the ideal 
$(F_1I,F_2,\ldots,F_{n-1})$ is called {\em basic double link}.

(1) Let $C\subset\mathbb P_K^4$ be the rational curve of degree $30$ 
parametrized by 
$$\left\{\begin{array}{ccc}
 x_0 &= & u^{30}+v^{30}\\
 x_1 &= & u^{29} v + u^{19} v^{11} + u^{9} v^{21} \\
 x_2 &= & u^{29} v + u^{18} v^{12} + u^{8} v^{22}\\
 x_3 &= & u^{27} v^{3} + u^{17} v^{13} + u^{7}v^{23}\\
 x_4 &= & u^{26} v^{4} + u^{16} v^{14} + u^{6}v^{24}\end{array}
\right..$$
Over a field of characteristic $31991$, we get $\rho_C+1=21=\omega_I$. 
Since $\rho_Z\leq\lceil \frac{deg(C)-1}{n-1} \rceil=10$, it is 
$reg(C)=21$ by Theorem 3.4(1). Further computations tell us that $\rho_Z=4$.

(2) Theorem 3.4(2) holds not only for aCM curves. Let $C_0\subset\mathbb P^3_K$ 
be a general elliptic curve of degree $5$. By applying a sequence of two basic 
double linkages to $C_0$, respectively of types $(1,5)$ and $(1,7)$, we obtain 
a non-aCM curve $X$ of degree $17$ with $reg(X)=7=\rho_{X}+2$, 
both in characteristics $31991$ and $0$. In fact, $\rho_Z=6$.

(3) Let $C_0\subset \mathbb P^3_K$ be a general rational curve of degree $4$.
By applying a sequence of two basic double linkages to $C_0$, respectively
of types $(1,4)$ and $(1,6)$, we obtain a curve $X$ of degree $14$ 
with $reg(X)=6=\rho_{X}+3$, both in characteristics $31991$ and $0$. In fact,
$\rho_Z=5$.
\end{example}

\begin{remark} \rm Let $C\subset \mathbb P^n_K$ be a connected non-special 
curve (i.e.
$h^1(\mathcal O_C(t))=0$ for every $t\geq 1$) over a field of characteristic
zero. We have:
\begin{itemize}
\item[(a)] If $2=\omega_I\geq \rho_C+1$, then $reg(C)\leq 3$; in particular,
$reg(C)=2$ if and only if $C$ is the rational normal curve.
\item[(b)] Otherwise, $reg(C)=\max\{\rho_C+1,\omega_I\}$. Namely, by
Theorem 2.1 of \cite{S3} it follows that $\rho_Z\leq \omega_I$.
Thus, if $\omega_I<\rho_C+1$, then $reg(C)=\rho_C+1$
by Theorem 3.4(1). If $\omega_I \geq \rho_C+1$, then
$h^1({\mathcal I}_C(\omega_I-1))=h^2({\mathcal I}_C(\omega_I-1))
=h^1({\mathcal O}_C(\omega_I-1))=0$ and $h^2({\mathcal I}_C(\omega_I-2))
=h^1({\mathcal O}_C(\omega_I-2))=0$ (see Remark 2.7), so $reg(C)=\omega_I$.
\end{itemize}
For example, from $(b)$, we get that every general rational curve $C\subset
\mathbb P^n_K$, with
$deg(C)>n$, being of maximal rank and non-special, has $reg(C)=\rho_C+1$
(see also \cite{O} and the references therein).
\end{remark}

\section{Regularity bounds for curves}

From now on, we suppose that the hyperplane section $Z$ of a curve 
$C\subset\mathbb P^n_K$ is obtained by a hyperplane which is general among 
those intersecting the curve properly. The following statement 
is based on results of \cite{DGO} about Hilbert function under liaisons
and the symmetry related to Gorenstein rings' Hilbert function.

\begin{proposition} Let $\alpha_1\leq\ldots\leq\alpha_{n-1}$ minimal degrees
for which there exists a c.i. of type
$(\alpha_1,\ldots,\alpha_{n-1})$ containing the curve $C$. Hence
\begin{itemize}
\item[(i)] $reg(C)\leq \max\{\rho_C+1,\alpha_1+\ldots +\alpha_{n-1}-(n-1)+1\}$;
\item[(ii)] if $deg(C)<\prod_{i=1}^{n-1}\alpha_i$, then $Z$ is not a c.i. 
of type $(\alpha_1,\ldots,\alpha_{n-1})$ and thus  
$reg(C)\leq \max\{\rho_C+1,\alpha_1+\ldots +\alpha_{n-1}-(n-1)\}$.
\end{itemize}
\end{proposition}

\begin{proof} 
We will see that the statement is a straightforward consequence of results of 
\cite{DGO} and of Theorem 3.4. 
First of all, recall that such a c.i. exists always 
(see Theorem 3.14 of Chapter VI of \cite{K}).

Next, let $\beta_1\leq\ldots\leq\beta_{n-1}$ be minimal degrees for which there 
exists a c.i. $W\subset\mathbb P^{n-1}_K$ of type 
$(\beta_1,\ldots,\beta_{n-1})$ containing the general hyperplane section 
$Z\subset \mathbb P^{n-1}_K$. Then, it is well known that the regularity 
of the Hilbert function of $W$ is $\beta_1+\ldots+\beta_{n-1}-(n-1)$ and
$reg(W)=\beta_1+\ldots+\beta_{n-1}-(n-1)+1$. Moreover, 
by Theorem 3(a) of \cite{DGO}, it follows that, for every integer $t$ such 
that $0\leq t\leq reg(W)-1$, we have $\Delta H_W(t)=\Delta H_(reg(W)-1-t)$, so 
in particular $\Delta H_W(reg(W)-1)=1$.
The ideal $(I(W):I(Z))$ is saturated (Lemma 5.2.1, \cite{Mi}) 
and defines a scheme $Y$ algebraically 
linked to $Z$ by $W$. By Theorem 3(b) of \cite{DGO} it is $\Delta H_W(t) 
= \Delta H_Z(t)+\Delta H_{Y}(reg(W)-1-t)$, for $0\leq t\leq reg(W)-1$, so 
in particular $\Delta H_W(reg(W)-1)= \Delta H_Z(reg(W)-1)+\Delta H_{Y}(0)$.
If $Z\neq W$, then $Y$ is not empty 
and hence $\Delta H_Y(0)=1$. So, $\Delta H_Z(reg(W)-1)=0$ and 
\begin{equation}
reg(Z) \leq reg(W)-1 = \beta_1+\ldots+\beta_{n-1}-(n-1).
\end{equation}
By the hypotheses, we have minimal degrees $\alpha_1,\ldots,\alpha_{n-1}$ for
which there exists a c.i. in $\mathbb P^n_K$ of type 
$(\alpha_1,\ldots,\alpha_{n-1})$ containing the curve $C$. Since regular 
sequences are preserved under general hyperplane sections, then we have
also a c.i. in $\mathbb P^{n-1}_K$ of type 
$(\alpha_1,\ldots,\alpha_{n-1})$ containing $Z$. The degrees of this complete 
intersection could be minimal or not for $Z$. Anyway we can take a complete 
intersection containing $Z$ of minimal degrees $\beta_1,\ldots,\beta_{n-1}$ 
such that $\beta_i\leq \alpha_i$, for every $1\leq i \leq n-1$. 
Hence, by Theorem 3.1, part (i) follows.

If $deg(Z)=deg(C)<\prod_{i=1}^{n-1}\alpha_i$, then $Z$ is not a c.i. 
of type $(\alpha_1,\ldots,\alpha_{n-1})$; thus $reg(Z)\leq 
\alpha_1+\ldots+\alpha_{n-1}-(n-1)$ and also part (ii) holds.
\end{proof}

\begin{example}\rm For every case of Theorem 3.4, we give examples of curves 
in $\mathbb P^4_K$ for which the bound of Proposition 4.1 is sharp. 
\begin{itemize} 
\item[(1)(a)] Let $C$ be a general degree $6$ elliptic curve, then $\rho_C=2$,
$reg(C)=3$, $\alpha_1=\alpha_2=\alpha_3=2$ and, so, 
$reg(C)=\alpha_1+\alpha_2+\alpha_3-(4-1)$.
\item[(b)] The ideal 
$I:=(x_0^2,x_1^3,x_2^4,x_1x_2^2x_3,x_0x_1^2x_2)\subset K[x_0,\ldots,x_4]$ 
is saturated and defines a curve $C\subset \mathbb P_K^4$ of degree $15$ 
such that $reg(C)=6=\rho_C+1=\rho_Z+2$ and 
$\alpha_1=2$, $\alpha_2=3$, $\alpha_3=4$.
\item[(2)(a)] Let $C$ be a general degree $5$ elliptic curve (which is aCM), 
then $\rho_C=1$, $reg(C)=3$, $\alpha_1=\alpha_2=\alpha_3=2$ and, so,
$reg(C)=\alpha_1+\alpha_2+\alpha_3-(4-1)$.
\item[(b)] The ideal
$I:=(x_0^2,x_1^3,x_2^4,x_0x_1,x_0x_2,x_0x_3)\subset K[x_0,\ldots,x_4]$
is saturated and defines a non-aCM curve $C\subset \mathbb P_K^4$
of degree $12$ such that $reg(C)=6=\rho_C+2=\rho_Z+1$ and 
$\alpha_1=2$, $\alpha_2=3$, $\alpha_3=4$. By further computation we note
that $Z$ is a degenerated c.i. of type $(3,4)$, 
hence $Z\subset \mathbb P^2_K$. 
\item[(3)\phantom{(a)}] The ideal
$I:=(x_0^2,x_1^3,x_2^3,x_0x_1,x_0x_3)\subset K[x_0,\ldots,x_4]$
is saturated and defines a curve $C\subset \mathbb P_K^4$
of degree $9$ such that $reg(C)=5=\rho_C+3=\rho_Z+1$ and
$\alpha_1=2$, $\alpha_2=3$, $\alpha_3=3$. As in (2), by further computation 
we note that $Z$ is a degenerated c.i. of type $(3,3)$,
hence $Z\subset \mathbb P^2_K$.
\end{itemize} 
\end{example}
When $C$ is a curve in $\mathbb P^3_K$ with hyperplane section 
$Z\subset \mathbb P^2_K$ we can have further information about 
$\rho_Z=reg(Z)-1$. Indeed, 
we want to improve and compare the bound of Proposition 4.1 with other bounds
of the same type for space curves. So, from now on we suppose that 
$C$ is a curve in $\mathbb P^3_K$ with general hyperplane section 
$Z\subset \mathbb P^2_K$. 
The following facts are straightforward refinements of Proposition 4.1 
for curves $C$ in $\mathbb P^3_K$.

\begin{corollary} 
{\rm (a)} If $C\subset \mathbb P^3_K$ is integral and $f_1,\ldots,f_s$ 
are minimal generators of $I$ with $d_1=deg(f_1)\leq \ldots\leq 
d_s=deg(f_s)$, then Proposition 4.1 holds with $\alpha_1=d_1$, $\alpha_2=d_2$.

{\rm (b)}
Let $C$ be equidimensional and locally Cohen-Macaulay over $K$ with 
characteristic $0$ and suppose $C$ to be a non-complete intersection 
with degree $> 4$. Then, if either
$deg(C)$ is odd or $C$ is not contained in a quadric, we have 
$reg(C) \leq \max\{\rho_C+1,\alpha_1+\alpha_2-2\}$. 
\end{corollary}

\begin{proof}
(a) Being $C$ integral, minimal generators of $I$ are always irreducible.
Hence $f_i, f_j$ form an $S$-regular sequence for every integers $i, j$ such
that $1\leq i < j\leq s$.

(b) By Theorem 2.3.1 of \cite{Mi} (that generalizes a result of R. Strano)
we have that a general hyperplane section $Z$ cannot be a complete
intersection. Then it is enough to apply Proposition 4.1(ii).
\end{proof}

To obtain interesting examples of integral curves, often we will apply 
known deformation techniques with constant cohomology. In these cases
we do not compute the curve explicitely and, so, we do not have further
information about the hyperplane section $Z$.

\begin{example}\rm
We give two examples for which the bound of Corollary 4.3(b) is sharp.
In the first example we exhibit a family of curves which are integral, 
while in the second example the given curve is not integral. 

(1)  Let $C_0$ be the curve of Moh \cite{Moh} parametrized (over a field
of characteristic $0$) by
$$\left\{\begin{array}{ccc}
  x_0 &= & u^{31}\\
  x_1 &= & u^{25} v^{6} + v^{31} \\
  x_2 &= & u^{23} v^{8}\\
  x_3 &= & u^{21} v^{10}\end{array}
\right..$$
By applying a sequence of two basic double linkages of type, respectively,
$(1,5)$ and $(w,7)$, with $w\geq 1$, we obtain a curve $X_w$ for which, by
Proposition 3.5 of \cite{No}, can be deformed with constant cohomology to an 
integral curve $C_w$ of degree $36+7w$ such that $\rho_{C_w}+1=10+w=reg(C_w)$.
Moreover, we have $\alpha_w=5+w$ and $\beta_w=7$ if $w=1,2$, whereas
$\alpha_w=7$ and $\beta_w=5+w$ otherwise. Thus,
$reg(C_w)=\max\{\rho_{C_w}+1,\alpha_w+\beta_w-2\} <
\max\{\rho_{C_w}+1,\lceil\frac{deg(C_w)-1}{2}\rceil+1\}$,
for every $w\geq 1$.

(2) Let $C_0\subset\mathbb P_K^3$ be the following smooth rational curve 
(over a field of characteristic $0$)
$$\left\{\begin{array}{ccc}
 x_0 &= & u^{12}+v^{12}\\
 x_1 &= & u^{11} v + u v^{11} \\
 x_2 &= & u^{10} v^{2} + v^{12}\\
 x_3 &= & u^{9} v^{3}\end{array}
\right.$$
to which we apply a basic double linkage of type $(1,9)$, obtaining a 
curve $C$ of degree $21$ such that 
$P_C(t)=21t-38$ and $reg(C)=12=\rho_{C}+1$, $\alpha_{I}=\alpha_1=5$, 
$\alpha_2=9$, $\omega_{I}=12$. So it happens that 
$\max\{\rho_{C}+1,\alpha_1+\alpha_2-2\}=12=reg(C)$. In this case we can also 
compute that $\rho_Z=8$ and that minimal degrees of a c.i. 
containing $Z$ are also $5$ and $9$.
\end{example}

For integral space curves over a field of characteristic $0$, 
one can say something more basing on the shape of Borel 
ideals. Namely, it is known that (see, for example, \cite{Gr}) the Borel ideal 
$gin(I(Z))$ is of the following type 
$$(x_0^s,x_0^{s-1}x_1^{\lambda_{s-1}},\ldots,x_0x_1^{\lambda_1},
x_1^{\lambda_0}),$$
where $\lambda_0 = \rho_Z+1$, $\lambda_i-2 \leq \lambda_{i+1}\leq\lambda_i-1$, 
for $0\leq i \leq s-1$, $s= \alpha_{I(Z)}$, 
$deg(Z)=\sum_{i=0}^{s-1} \lambda_i$ and
(Corollary 4.9 of \cite{Gr})
\begin{equation}
\frac{deg(Z)}{s}+\frac{s-1}{2}\leq \lambda_0\leq
\frac{deg(Z)}{s}+s-1.
\end{equation}

\begin{proposition} Let $\gamma_C:=\max\{t\in \mathbb N : deg(C)>t^2+1\}$; then
\begin{itemize}
\item[(a)] $deg(C)>\alpha_I^2+1$ $\Rightarrow$ $reg(C)\leq 
\max\{\rho_C+1,\lfloor\frac{deg(C)}{\alpha_I}\rfloor+\alpha_I-1\}$.
\item[(b)] $deg(C)\leq\alpha_I^2+1$ $\Rightarrow$ $reg(C)\leq
\max\{\rho_C+1,\lfloor\frac{deg(C)}{\gamma_C+1}\rfloor+\gamma_C\}$.
\end{itemize}
\end{proposition}

\begin{proof}
{\it (a)} By Theorem 3.4 the thesis is an easy consequence of formula (4.2) 
and Laudal's lemma of \cite{S2}, because one obtains that 
$s=\alpha_{\bar J^{sat}}=\alpha_J = \alpha _I$.

{\it (b)} In this case $s$ is not known. However we obtain that $s>\gamma_C$,
namely, in our hypothesis it can happen that 
either $s=\alpha_I>\gamma_C$ or $2\leq s <\alpha_I$. In the second case, 
by \cite{S2} $deg(C)\leq s^2+1$ and so $s>\gamma_C$.
Therefore $deg(C) \geq \lambda_0+\lambda_1+\ldots+\lambda_{\gamma_C} \geq
\lambda_0+\lambda_0-2+\ldots+\lambda_0-2\gamma_C =
(\gamma_C+1)(\lambda_0-\gamma_C)$ and we are done.
\end{proof}

\begin{example} \rm We give examples for which the described bounds are sharp.

(1) By Proposition 3.5 of \cite{No}, the curve $X$ of Example 3.6(2) can be 
deformed with constant cohomology to an integral curve $C$ of degree $17$, 
regularity $reg(C)=7=\rho_{C}+2$, initial degree 
$\alpha_I=\alpha_1=5=\alpha_2$. Thus 
$\gamma_C=3$ and $reg(C)=7=\max\{\rho_{C}+1,\lfloor\frac{deg(C)}{4}\rfloor+3\}$, 
whereas $\max\{\rho_{C}+1,\alpha_1+\alpha_2-2\}=8$ and
$\max\{\rho_{C}+1,\lceil\frac{deg(C)-1}{2}\rceil+1\}=9$.

(2) By Proposition 3.5 of \cite{No}, the curve $X$ of Example 3.6(3) can be
deformed with constant cohomology to an integral curve $C$ of degree $14$, 
regularity $reg(C)=6=\rho_{C}+3$, initial degree $\alpha_{I}=\alpha_1=4$ and 
$\alpha_2=5$. Thus $\gamma_C=3$ and $reg(C)=6=
\max\{\rho_{C}+1,\lfloor\frac{deg(C)}{4}\rfloor+3\}$, whereas 
$\max\{\rho_{C}+1,\alpha_1+\alpha_2-2\}=7$ and $\max\{\rho_{C}+1,
\lceil\frac{deg(C)-1}{2}\rceil+1\}=8$.

(3) By applying a further basic double linkage of type $(1,4)$ to the curve $C$ 
in (2) above, we obtain a curve $X'$ which, by Proposition 3.5 of \cite{No}, 
can be deformed with constant cohomology to an integral curve $C'$ (defined 
by an ideal $I'$) of degree $18$ with regularity $reg(C')=7=\rho_{C'}+3$,
initial degree $\alpha_{I'}=\alpha_1=4$ and $\alpha_2=6$. Thus $reg(C')=7= 
\max\{\rho_{C'}+1,\lfloor\frac{deg(C')}{\alpha_{I'}}\rfloor+\alpha_{I'}-1\}$, 
whereas $\max\{\rho_{C'}+1,\alpha_1+\alpha_2-2\}=8$ and
$\max\{\rho_{C'},\lceil\frac{deg(C')-1}{2}\rceil+1\}=10$.

(4) In Example 3.2 of \cite{MM-R}, over a field of characteristic $0$,
a family $\{Y_m\}_{m\geq 0}$ of $2$-Buchsbaum smooth integral curves is 
constructed such that $reg(Y_m)=\max\{m+4,2m+3\}$, the initial degree is 
$\alpha_m=\alpha_1=\alpha_2=m+3$, $deg(Y_m)=m^2+4m+6 < \alpha_m^2+1$ and 
$\gamma_m:=\gamma_{Y_m}=m+2$. 

If $m=0$, $Y_0$ is a general elliptic curve of degree $6$ and 
$4=reg(Y_0)= \rho_{Y_0} +1$.

If $m=1$, then $5=\rho_{Y_1}+1=reg(Y_1) = 
\max\{\rho_{Y_1}+1,\lfloor\frac{deg(Y_1)}{4}\rfloor+3\}<6 = 
\max\{\rho_{Y_1}+1,\alpha_1+\alpha_2-2\}<8=
\max\{\rho_{Y_1}+1,\lceil\frac{deg(Y_1)-1}{2}\rceil+1\}$.

If $m\geq 2$, then $\rho_{Y_m}+2=reg(Y_m)=2m+3 
=\max\{\rho_{Y_m}+1,\lfloor\frac{deg(Y_m)}{m+3}\rfloor+m+2\}<
2m+4=\max\{\rho_{Y_m}+1,\alpha_1+\alpha_2-2\}<
\max\{\rho_{Y_m}+1,\lceil\frac{deg(Y_m)-1}{2}\rceil+1\}$.
\end{example}

\medskip

We thank U. Nagel whose suggestions inspired this paper.

\providecommand{\bysame}{\leavevmode\hbox to3em{\hrulefill}\thinspace}
\providecommand{\MR}{\relax\ifhmode\unskip\space\fi MR }
\providecommand{\MRhref}[2]{%
  \href{http://www.ams.org/mathscinet-getitem?mr=#1}{#2}
}
\providecommand{\href}[2]{#2}


\end{document}